\documentclass[12pt]{article}
\usepackage{amsmath,amscd,amsthm,amssymb}

\usepackage{color}
\newcommand{\tcr}{\textcolor{red}}

\newtheorem{theorem}{Theorem}{}
\newtheorem{corollary}{Corollary}{}
\newtheorem{definition}{Definition}{}
{}
\newtheorem{lemma}{Lemma}{}
{}
\newcommand{\es}{\mathop{\rm ess \; inf}\limits}
\newcommand{\ess}{\mathop{\rm ess \; sup}\limits}
\def\R{\mathbb{R}}
\def\Rn{{\mathbb{R}^n}}
\def\i{\infty}
\def\Lloc{L_1^{\rm loc}(\Rn)}

\chardef\No=24

\begin{document}

\begin{center}
\Large \bf   Maximal and Calder\'{o}n-Zygmund operators on the local variable Morrey-Lorentz spaces and some applications
\end{center}

\centerline{\large A. Kucukaslan$^{a,b,*,}$\footnote{
The research of A. Kucukaslan was supported by the grant of The Scientific and Technological Research
Council of Turkey (2219 TUBITAK, Grant-1059B191600675).
\\
The research of V.S. Guliyev was partially supported by the grant of 1st Azerbaijan-Russia Joint Grant Competition (Agreement Number No. EIF-BGM-4-RFTF-1/2017-21/01/1-M-08).
\\
The research of V.S. Guliyev and A. Serbetci was partially supported by the grant of Cooperation Program 2532 TUBITAK-RFBR (Russian foundation for basic research) with Agreement Number No. 119N455.
\\
{$*$ Corresponding author.}
\\
E-mail addresses: kucukaslan@pau.edu.tr (A. Kucukaslan), vagif@guliyev.com (V.S. Guliyev), aykol@science.ankara.edu.tr (C. Aykol), serbetci@ankara.edu.tr (A. Serbetci).
}, V.S. Guliyev$^{c,d,e}$, C. Aykol$^{f}$, A. Serbetci$^{f}$}

\

\centerline{$^{a}$\it Institute of Mathematics, Czech Academy of Sciences, Prague, Czech Republic}
\centerline{$^{b}$\it Pamukkale University, School of Applied Sciences, Denizli, Turkey}
\centerline{$^{c}$\it Institute of Applied Mathematics, Baku State University, Baku, Azerbaijan}
\centerline{$^{d}$\it Institute of Mathematics and Mechanics, Baku, Azerbaijan}
\centerline{$^{e}$\it Dumlupinar University, Department of Mathematics, Kutahya, Turkey}
\centerline{$^{f}$\it Ankara University, Department of Mathematics, Ankara, Turkey}
\

\begin{abstract} In this paper, we give the definition of local variable Morrey-Lorentz spaces $\mathcal{M}_{p(\cdot),q(\cdot),\lambda}^{loc}(\Rn)$ which are a new class of functions. Also, we prove the boundedness of the Hardy-Littlewood maximal operator $M$ and Calder\'{o}n-Zygmund operators $T$ on these spaces including the class of sublinear operators $T_0$ generated by Calder\'{o}n-Zygmund operators.
Finally, we apply these results to the Bochner-Riesz operator $B_r^\delta$, identity approximation $A_{\varepsilon}$ and the Marcinkiewicz operator $\mu_\Omega$ on the spaces $\mathcal{M}_{p(\cdot),q(\cdot),\lambda}^{loc}(\Rn)$.

\end{abstract}

\



\

\section{Introduction}

The study of function spaces with variable exponent has been stimulated by problems of elasticity, fluid dynamics, calculus of variations and differential equations with non-standard growth conditions (see \cite{RX, DHHR, Ru, RY,  Zh}). Various results on non-weighted and weighted boundedness in variable exponent Lebesgue spaces have been proved for maximal, singular and fractional type operators, we refer to surveying papers \cite{DHN} and \cite{Samko}.
In \cite{EKS} variable exponent Lorentz spaces $L_{p(\cdot),q(\cdot)}$ are introduced and the boundedness of the singular integral and fractional type operators and corresponding ergodic operators are proved in these spaces.

The Lorentz-Morrey space $\mathcal{L}_{p,q;\lambda}(\Rn)$ was first defined in \cite{Mg} and also considered in \cite{Hat, Ho, Ragusa2012}. Later, the local Morrey-Lorentz spaces $\mathcal{M}_{p,q;\lambda}^{loc}(\Rn)$ are introduced and the basic properties of these spaces are given in \cite{AGS}. These spaces are a very natural generalization of the Lorentz spaces such that $\mathcal{M}_{p,q;0}^{loc}(\Rn)=L_{p,q}(\Rn)$. Recently, in \cite{AGKS, GAKS} and \cite{GKAS} the authors have studied the boundedness of the Hilbert transform, the Hardy-Littlewood maximal operator $M$ and the Calder\'{o}n-Zygmund operators $T$, and the Riesz potential $I_\alpha$ on the local Morrey-Lorentz spaces $\mathcal{M}_{p,q;\lambda}^{loc}$ by using related rearrangement inequalities, respectively. In \cite{Miz}, the authors give  the definition of central Lorentz-Morrey space of variable exponent by the symmetric decreasing rearrangement. They prove the boundedness of maximal operator in these spaces and establish Sobolev's inequality for Riesz potentials.

In this present paper, we define the local variable Morrey-Lorentz spaces 
$\mathcal{M}_{p(\cdot),q(\cdot),\lambda}^{loc}(\Rn)$ and prove the boundedness of the Hardy-Littlewood maximal operator $M$ and Calder\'{o}n-Zygmund operators $T$ on these spaces including the class of sublinear operators $T_0$ generated by Calder\'{o}n-Zygmund operators. We also give some applications of our results.

The paper is organized as follows. In Section 2, we give some notation and definitions.
We introduce local variable Morrey-Lorentz spaces $\mathcal{M}_{p(\cdot),q(\cdot);\lambda}^{loc}(\Rn)$.
In Section 3, we prove the boundedness of the maximal operator and the Calder\'{o}n- Zygmund operators $T$ in the spaces $\mathcal{M}_{p(\cdot),q(\cdot),\lambda}^{loc}(\Rn)$. In Section 4,
finally, as applications we get the boundedness of Bochner-Riesz operator $B_r^\delta$, Marcinkiewicz operator $\mu_\Omega$ and the sublinear operator $\sup_{\varepsilon>0}|A_{\varepsilon}f(x)|$ on the spaces $\mathcal{M}_{p(\cdot),q(\cdot),\lambda}^{loc}(\Rn)$, where $A_{\varepsilon}$ is the identity approximation.

Throughout the paper we use the letter $C$ for a positive constant, independent of appropriate parameters and not necessary the same at each occurrence.

\section{Preliminaries}
For $x \in \Rn $ and $r > 0,$ let  $B(x,r)$ denote the open ball centered at $x$ of radius $r$ and $|B(x,r)|$ be the Lebesgue measure of the ball $B(x,r)$.
Note that, $|B(x,r)|=\omega_{n}r^{n}$, where $\omega_{n}$ is the volume of the unit ball in $\Rn$.

Let $f$ be a locally integrable function on $\Rn$. Hardy-Littlewood maximal function $Mf$ is defined by
\begin{equation*}
Mf(x)=\sup_{r>0}\frac{1}{|B(x,r)|}\int_{B(x,r)}|f(y)|dy, ~~x \in \Rn.
\end{equation*}
Maximal operators play an important role in the differentiability properties of functions, singular integrals and partial differential equations. They often provide a deeper and more simplified approach to understand problems in these areas. For the operator $M$ the rearrangement inequality
\begin{equation}\label{cay0}
cf^{\ast\ast}(t) \leq (Mf)^{\ast}(t)\leq Cf^{\ast\ast}(t), ~~t \in (0,\infty)
\end{equation}
holds, where $c$ and $C$ are independent of $f$ and $f^{*}$ denotes the right continuous non-increasing rearrangement
of $f$:
\begin{equation*}
f^{*}(t):=inf \left\{{\lambda > 0: \mu_{f}(\lambda) \le t}\right\}, ~~   t\in (0,\infty ),
\end{equation*}
and \begin{equation*}
\mu_{f}(\lambda):=\left| \{y\in \R : |f(y)|>\lambda\}\right|
\end{equation*} is the distribution function of the function $f$.

Let $T$ be a Calder\'{o}n-Zygmund operator, i.e., a linear operator bounded from
$L_2(\Rn)$ to $L_2(\Rn)$ taking all infinitely continuously
differentiable functions $f$ with compact support to the functions
$T f \in \Lloc$ represented by
\begin{equation*} 
Tf(x) =  \int_{\Rn }K(x,y) f(y)dy,   ~~ x \notin{\rm supp} f,
\end{equation*}
provided it exists almost everywhere. Here $K(x,y)$ is a continuous function away from the diagonal which
satisfies the standard estimates: there exist $C>0$ and
$0<\varepsilon\leq 1$ such that
\begin{equation*} 
|K(x,y)| \le C |x-y|^{-n}
\end{equation*}
for all $x, y \in \Rn, \,\, x\neq y$, and
\begin{equation*} \label{SIOker1}
|K(x,y)-K(x',y)|+|K(y,x)-K(y,x')|\le C
\left(\frac{|x-x'|}{|x-y|}\right)^{\varepsilon}\, |x-y|^{-n} ,
\end{equation*}
whenever $2|x-x'|\le |x-y|$. More information about such operators can be found in
\cite{CM, Dynkin}.

The following rearrangement inequality
\begin{equation} \label{BS1}
(Tf)^{*}(t)\leq C\left(\int_{0}^{t} f^{*}(s) ds + \int_{t}^{\i} f^{*}(s) \frac{ds}{s}\right)
\end{equation}
is valid for the Calder\'{o}n-Zygmund operators $T$, where where $C$ is independent of $T$  (see e.g. \cite{BenRud}). Similar sharp rearrangement estimates are of great importance in the study of operators on rearrangement-invariant function spaces as well as in interpolation theory.

Suppose that $T_0$ represents a linear or a sublinear operator, such that that for any $f\in L_1(\Rn)$ with compact support and $x\notin supp f$
\begin{equation}\label{sublR}
	|T_0 f(x)|\le c_0 \int\limits_{\Rn} \frac{|\Omega(x-y)|}{|x-y|^{n}}~ |f(y)| dy,
\end{equation}
where $c_0$ is independent of $f$ and $x$.

We point out that the condition \eqref{sublR} was first introduced by Soria and Weiss in \cite{SW}.
The condition \eqref{sublR} are satisfied by many interesting operators in harmonic analysis, such as the Calder\'{o}n-Zygmund operators,
Carleson's maximal operator, Hardy--Littlewood maximal operator, C. Fefferman's singular multipliers, R. Fefferman's singular integrals, Ricci-Stein's oscillatory singular integrals, the Bochner--Riesz means and so on (see  \cite{GulAlKarSh,LLY,SW} for details).

Let $p(t)$ be a measurable function on $(0,\i)$. We suppose that
\begin{equation*}
1< p_-\le p(t)\le p_+<\infty,
\end{equation*}
where
\begin{equation*}
p_{-} := \es_{0<t<\i} \, p(t),  ~~  p_{+} := \ess_{0<t<\i} \, p(t).
\end{equation*}
We denote by $p'(\cdot)=\frac{p(t)}{p(t)-1}.$
We will use the following decay conditions:
\begin{equation} \label{cny10}
|p(t)-p(0)|\leq \frac{A_{0}}{|ln\,t|},\,\,\, 0<t\leq\frac{1}{2},
\end{equation}
\begin{equation} \label{cny11}
 |p(t)-p(\infty)|\leq \frac{A_{\infty}}{ln\,t},\,\,\, t\geq 2,
\end{equation}
where $A_0, A_{\i}>0$ do not depend on $t$.

By $p \in \mathcal{P}_{0,\i}(0,\i)$ we denote the set of bounded measurable functions (not necessarily with values in $[1,\i)),$ which satisfy the decay conditions $\eqref{cny10}$ and $\eqref{cny11}$. Also, by $L_{p(\cdot)}(0,\i)$ we denote the variable exponent Lebesgue space of measurable functions $\varphi$ on $(0,\i)$ such that
\begin{equation*}
 \mathcal{J}_{p(\cdot)}(\varphi)= \int_{0}^{\i}|\varphi(s)|^{p(s)}ds < \infty.
\end{equation*}
This is a Banach function space with respect to the norm (see e.g. \cite{ELN})
\begin{equation*}
\|\varphi\|_{L_{p(\cdot)}}=\inf \left\{{\lambda > 0: \mathcal{J}_{p(\cdot)} \left(\frac{\varphi}{\lambda}\right) \leq 1}\right\}.
\end{equation*}

\begin{definition}\cite{LPSW}\label{Ca}
Let  $1 < q_{-} \leq q_{+}< \infty$, $0 \leq \lambda < 1$. We denote by $LM_{q(\cdot),\lambda} \equiv LM_{q(\cdot),\lambda}(0,\infty)$
the variable local Morrey space with finite norm
\begin{align*}
\|\varphi\|_{LM_{q(\cdot),\lambda}} &=\sup_{r>0}r^{-\frac{\lambda}{q_{*}(r)}}\| \varphi\|_{L_{q(\cdot)}(0,r)}
\\
&=\sup_{r>0}\inf \left\{\eta > 0:\int_{0}^{r}\left|\frac{\varphi(s)}{\eta r^{\frac{\lambda}{q_{*}(r)}}}\right|^{q(s)}ds \leq 1 \right\},
\end{align*}
where $q_{*}(r)=q(0), 0<r<1$ and $q_{*}(r)=q(\infty), r\geq1$.
\end{definition}

\begin{definition}\cite{EKS}
 Let $1 \leq p_{-} \leq p_{+}< \infty,~ 1 < q_{-} \leq q_{+}< \infty$. We denote by  $L_{p(\cdot),q(\cdot)}(\Rn)$  variable Lorentz space, the space of functions $f$ on $\Rn$ such that $t^{\frac{1}{p(t)}-\frac{1}{q(t)}}f^{*}(t) \in L_{q(\cdot)}(0,\i) $, i.e.
\begin{equation*}\label{cay}
\mathcal{J}_{p(\cdot),q(\cdot)}(f)=\int_{0}^{\i}t^{\frac{q(t)}{p(t)}-1}(f^{*}(t))^{q(t)}dt< \i
\end{equation*}
and we denote
\begin{equation*}
\|f\|_{L_{p(\cdot),q(\cdot)}(\Rn)}=\inf \left\{{\sigma > 0: \mathcal{J}_{p(\cdot),q(\cdot)} \left(\frac{f}{\sigma}\right) \leq 1}\right\}=\left\| t^{\frac{1}{p(t)}-\frac{1}{q(t)}}f^{*}(t) \right\|_{L_{q(\cdot)}(0,\i)},
\end{equation*}
where $f^{*}$ denotes the non-increasing rearrangement
of $f$ such that
\begin{equation*}
f^{*}(t)=\inf \left\{{\lambda > 0: \left| \{y\in \Rn : |f(y)|>\lambda\}\right| \le t}\right\}, ~~  \forall t\in (0,\infty )
\end{equation*}
 and
\begin{equation*}
f^{**}(t)=\frac{1}{t}\int_0^t f^{\ast}(s)ds
\end{equation*}
 (see \cite{BenSh}).
More information about variable Lorentz spaces can be found in \cite{EKS,KV}.
\end{definition}
We find it convenient to define the local variable Morrey-Lorentz spaces in the form as following.
\begin{definition}\label{kuc01}
Let $0 < p_{-} \leq p_{+}< \infty ,~ 1 < q_{-} \leq q_{+}< \infty$ and $0 \leq \lambda < 1$. We denote by $\mathcal{M}_{p(\cdot),q(\cdot),\lambda}^{loc}\equiv \mathcal{M}_{p(\cdot),q(\cdot),\lambda}^{loc}(\Rn)$
the local variable Morrey-Lorentz space, the space of all measurable functions with finite quasinorm
\begin{align*}
\|f\|_{\mathcal{M}_{p(\cdot),q(\cdot),\lambda}^{loc}} &:=\sup_{r>0}r^{-\frac{\lambda}{q_{*}(r)}}\| t^{\frac{1}{p(t)}-\frac{1}{q(t)}}f^{*}(t)\|_{L_{q(\cdot)}(0,r)}.
\end{align*}

These spaces generalize variable Lorentz spaces such that $\mathcal{M}_{p(\cdot),q(\cdot);0}^{loc}=L_{p(\cdot),q(\cdot)}$, when $\lambda=0$ (see \cite{EKS}). Also, if $\lambda=0$ and $q(\cdot)=p(\cdot)$ then $\mathcal{M}_{p(\cdot),p(\cdot);0}^{loc}=L_{p(\cdot)}$ are variable Lebesgue spaces (see \cite{KS}).

\end{definition}

\section{The maximal operator $M$ and the Calder\'{o}n-Zygmund operators $T$ in the spaces $\mathcal{M}_{p(\cdot),q(\cdot),\lambda}^{loc}(\Rn)$}

In this section, we prove the boundedness of the maximal operator $M$ and the Calder\'{o}n-Zygmund operators $T$ in the local variable Morrey-Lorentz spaces.
We need the following two definitions about Hardy operators which are used in the
proof of our main theorems. These operators are very important in analysis and have been widely
studied.
\begin{definition}\cite{EKS}
Let $\varphi$ be a measurable function on $(0,\infty)$. The weighted Hardy operators $H_{\beta(\cdot)}$ and $\mathcal{H}_{\beta(\cdot)}$ with power weight acting on $\varphi$ are defined by
\begin{equation*}
H_{\beta(\cdot)}\varphi(t)=t^{\beta(t)-1}\int_{0}^{t}\frac{\varphi(s)}{s^{\beta(s)}}ds
\end{equation*}
and
\begin{equation*}
\mathcal{H}_{\beta(\cdot)}\varphi(t)=t^{\beta(t)}\int_{t}^{\infty}\frac{\varphi(s)}{s^{\beta(s)+1}}ds.
\end{equation*}
\end{definition}
The following lemma provides some minimal assumptions on the function $r^{\frac{\lambda}{q_{*}(r)}}$ under which the so-defined spaces contain "nice" functions.

\begin{lemma}\label{kuc21} \cite{LPSW}
Let $1 < q_{-}\leq q_{+}<\infty$, $q\in \mathcal{P}_{0}(0,\i)$, $[r]_{1}=min\{1,r\}$ and $0 \leq \lambda < 1$. Then the assumption
\begin{equation}\label{kuc22}
\sup_{r>0} \, [r]_{1}^{\frac{1}{q(0)}} \, r^{-\frac{\lambda}{q_{*}(r)}}<\i
\end{equation}
is sufficient for bounded functions $f$ with compact support  to belong to the local variable Morrey spaces $LM_{q(\cdot),\lambda}(\Rn)$.
\end{lemma}

\begin{lemma}\label{cay0001} \cite{LPSW} Let $1 < q_{-}\leq q_{+}<\infty$,  $q\in \mathcal{P}_{0,\i}(0,\i)$, $0 \leq \lambda < 1$,  $\lim_{t\rightarrow 0}t^{\beta(t)}$ exists and finite and the condition \eqref{kuc22} satisfies.
 Suppose that the following conditions hold.
\


$(i)$ $t^{\beta(t)-a}$ and  $t^{\frac{\lambda}{q_{*}(t)}-\beta(t)-a}$ are almost decreasing for some $a\in \R$, in the case of operator $H_{\beta}$.
\

\

$(ii)$ $t^{-\beta(t)+b}$ and $t^{\frac{\lambda}{q_{*}(t)}-\beta(t)+b}$ are  almost increasing for some $b\in \R$, in the case of operator $\mathcal{H}_{\beta}$.

Then the conditions
\begin{equation*}\label{cay001}
\beta(t)<\frac{\lambda}{q_{*}(t)}+\frac{1}{q'(0)},~~~~~~ \beta(t)>\frac{\lambda}{q_{*}(t)}-\frac{1}{q(\i)}
\end{equation*}
are sufficient for the Hardy operators $H_{\beta}$ and $\mathcal{H}_{\beta}$, respectively, to be defined on the space $LM_{q(\cdot),\lambda}(\Rn)$.
\end{lemma}

\begin{lemma} \label{cay002}\cite{LPSW} Let  $1 < q_{-}\leq q_{+}<\infty$,  $q\in \mathcal{P}_{0,\i}(0,\i)$, $0 \leq \lambda < 1$. Suppose also that the conditions \eqref{kuc22} and of Lemma \ref{cay0001} are satisfied. Then the operators $H_{\beta}$ and $\mathcal{H}_{\beta}$ are bounded in the space $LM_{q(\cdot),\lambda}(0, \i)$ if $\beta(t)<\frac{\lambda}{q_{*}(t)}+\frac{1}{q'(t)}$, $\beta(t)>\frac{\lambda}{q_{*}(t)}-\frac{1}{q(t)}$, respectively.
\end{lemma}

The following theorem is one of the main results of our paper in which we give the boundedness of maximal operator in the local variable Morrey-Lorentz spaces.

\begin{theorem}\label{cay4} Let   $1\leq p_{-}\leq p_{+}<\infty$, $1 < q_{-}\leq q_{+}<\infty$, $p,q\in \mathcal{P}_{0,\i}(0,\i)$,  $0 \leq \lambda < 1$ and $f\in \mathcal{M}_{p(\cdot),q(\cdot),\lambda}^{loc}(\Rn)$. Suppose that the conditions \eqref{kuc22} and of Lemma \ref{cay0001} are satisfied. Then the maximal operator $M$ is bounded on the local variable Morrey-Lorentz spaces $\mathcal{M}_{p(\cdot),q(\cdot),\lambda}^{loc}(\Rn)$.
\end{theorem}

\begin{proof}Let $1\leq p_{-}\leq p_{+}<\infty$, $1 < q_{-} \leq q_{+}< \infty $, $0 \leq \lambda < 1$, the condition \eqref{kuc22} satisfies and $f\in \mathcal{M}_{p(\cdot),q(\cdot),\lambda}^{loc}(\Rn)$.  From the definition of local variable Morrey-Lorentz spaces and the inequality \eqref{cay0} we get
\begin{align*}\|Mf\|_{\mathcal{M}_{p(\cdot),q(\cdot),\lambda}^{loc}(\Rn)}
&=\sup_{r>0}r^{-\frac{\lambda}{q_{*}(r)}}\left\|t^{\frac{1}{p(t)}-\frac{1}{q(t)}}(Mf)^{*}(t)\right\|_{L_{q(\cdot)}(0,r)}
\\
&\leq C\sup_{r>0}r^{-\frac{\lambda}{q_{*}(r)}}\left\|t^{\frac{1}{p(t)}-\frac{1}{q(t)}}f^{**}(t)\right\|_{L_{q(\cdot)}(0,r)}
\\
&=C\sup_{r>0}r^{-\frac{\lambda}{q_{*}(r)}} \Big\|t^{\frac{1}{p(t)}-\frac{1}{q(t)}-1}\int_{0}^{t}f^{*}(s)ds\Big\|_{L_{q(\cdot)}(0,r)}
\\
&=C\|H_{\beta}g \|_{LM_{q(\cdot),\lambda}(0,\infty)},
\end{align*}
where $g(t)=t^{{\frac{1}{p(t)}-\frac{1}{q(t)}}}f^{*}(t), \beta(t)={\frac{1}{p(t)}-\frac{1}{q(t)}}$.

For $\beta(t)={\frac{1}{p(t)}-\frac{1}{q(t)}}$
the inequality $\beta(t)<\frac{1}{q'(t)}+\frac{\lambda}{q_{*}(t)}$ holds.
Therefore by Lemma \ref{cay002} we get
\begin{align}\label{cay2}
\|H_{\beta} g\|_{LM_{q(\cdot),\lambda}(0,\infty)}
&\leq C \|g\|_{LM_{q(\cdot),\lambda}(0,\infty)}\notag
\\
&=C\sup_{r>0}r^{-\frac{\lambda}{q_{*}(r)}}
\|t^{\frac{1}{p(t)}-\frac{1}{q(t)}}f^{*}(t)\|_{L_{q(\cdot)}(0,r)}\notag
\\
&=C\|f\|_{\mathcal{M}_{p(\cdot),q(\cdot),\lambda}^{loc}(\Rn)}.
\end{align}
From \eqref{cay2} we obtain the boundedness of the maximal operator $M$ in the space ${\mathcal{M}_{p(\cdot),q(\cdot),\lambda}^{loc}}(\Rn)$.
\end{proof}

In the case $\lambda=0,$ from Theorem \ref{cay4} we get the boundedness of the maximal operator $M$ in the variable Lorentz spaces $L_{p(\cdot),q(\cdot)}(\Rn)$ which is proved in \cite{EKS}.




In the following theorem is the other main result of our paper in which we give the boundedness of Calder\'{o}n-Zygmund operators in the local variable Morrey-Lorentz spaces.

\begin{theorem}\label{cay5} Let   $1\leq p_{-}\leq p_{+}<\frac{q_*(t)}{\lambda}$, $1\leq q_{-}\leq q_{+}<\infty$, $p,q\in \mathcal{P}_{0,\i}(0,\i)$, $0 \leq \lambda < 1$ and $f\in \mathcal{M}_{p(\cdot),q(\cdot),\lambda}^{loc}(\Rn)$. Suppose that the conditions \eqref{kuc22} and of Lemma \ref{cay0001} are satisfied. Then the Calder\'{o}n-Zygmund operator $T$ exists
almost everywhere $x \in \Rn$. Moreover, $T$ is  bounded on the local variable Morrey-Lorentz spaces $\mathcal{M}_{p(\cdot),q(\cdot),\lambda}^{loc}(\Rn)$.
\end{theorem}

\begin{proof}
Let $1\leq p_{-}\leq p_{+}<\frac{q_*(t)}{\lambda}, 1 \leq q_{-} \leq q_{+}< \infty $, $0 \leq \lambda < 1$, the condition \eqref{kuc22} satisfies and $f\in \mathcal{M}_{p(\cdot),q(\cdot),\lambda}^{loc}(\Rn)$.
From the definition of norm in local variable Morrey-Lorentz spaces and by using the inequality \eqref{BS1} we get
\begin{align*}
\|Tf\|_{\mathcal{M}_{p(\cdot),q(\cdot),\lambda}^{loc}(\Rn)}&=
\sup_{r>0}r^{-\frac{\lambda}{q_{*}(r)}}\|t^{\frac{1}{p(t)}-\frac{1}{q(t)}}(Tf)^{*}(t)\|_{L_{q(\cdot)}(0,r)}
\\
&\leq C\sup_{r>0}r^{-\frac{\lambda}{q_{*}(r)}} \Big\|t^{\frac{1}{p(t)}-\frac{1}{q(t)}-1}\int_{0}^{t}f^{*}(s)ds \Big\|_{L_{q(\cdot)}(0,r)}
\\
&~~~~~~~~+C\sup_{r>0}r^{-\frac{\lambda}{q_{*}(r)}} \Big\|t^{\frac{1}{p(t)}-\frac{1}{q(t)}}\int_{t}^{\i}\frac{f^{*}(s)}{s}ds \Big\|_{L_{q(\cdot)}(0,r)}
\\
&=I_{1}+I_{2}.
\end{align*}
$I_{1}$ can be estimated using the same method as in the proof of the boundedness of the maximal operator on $\mathcal{M}_{p(\cdot),q(\cdot),\lambda}^{loc}(\Rn)$ in Theorem \ref{cay4}.

Let us estimate $I_{2}:$
\begin{align*}
I_{2}&= C\sup_{r>0}r^{-\frac{\lambda}{q_{*}(r)}} \Big\|t^{\frac{1}{p(t)}-\frac{1}{q(t)}}\int_{t}^{\i}\frac{f^{*}(s)}{s}ds \Big\|_{L_{q(\cdot)}(0,r)}\notag
\\
&=C\|\mathcal{H}_{(\frac{1}{p(t)}-\frac{1}{q(t)})}g\|_{LM_{q(\cdot),\lambda}(0, \i)},
\end{align*}
where $g(t)=t^{{\frac{1}{p(t)}-\frac{1}{q(t)}}}f^{*}(t)$.
Since $\frac{1}{p(t)}-\frac{\lambda}{q_{*}(t)}>0$, for $\beta(t)={\frac{1}{p(t)}-\frac{1}{q(t)}}$
the inequality $\beta(t)>\frac{\lambda}{q_{*}(t)}-\frac{1}{q(t)}$ holds.
By Lemma \ref{cay002} we get
\begin{align}\label{cal2}
\|\mathcal{H}_{(\frac{1}{p(t)}-\frac{1}{q(t)})}g\|_{LM_{q(\cdot),\lambda}(0, \i)} &\leq C\|g\|_{LM_{q(\cdot),\lambda}(0, \i)}\notag
\\
&=C\sup_{r>0}r^{-\frac{\lambda}{q_{*}(r)}} \|t^{\frac{1}{p(t)}-\frac{1}{q(t)}}f^{*}(t)\|_{L_{q(\cdot)}(0,r)}\notag
\\
&=C\|f\|_{\mathcal{M}_{p(\cdot),q(\cdot),\lambda}^{loc}(\Rn)}.
\end{align}
Therefore we get $I_{2} \leq C \|f\|_{M_{p,q;\lambda}^{loc}(\Rn)}$.
Consequently we obtain the boundedness of $T$ in $\mathcal{M}_{p(\cdot),q(\cdot),\lambda}^{loc}(\Rn)$ from the inequalities \eqref{cay2} and \eqref{cal2} .
\end{proof}

In the case $\lambda=0,$ from Theorem \ref{cay5} we get the boundedness of the Calder\'{o}n-Zygmund operator $T$ in the variable Lorentz spaces $L_{p(\cdot),q(\cdot)}(\Rn)$ which is proved in \cite{EKS}.


In the following theorem we give the boundedness of sublinear operators $T_0$ generated by Calder\'{o}n-Zygmund operators in  $\mathcal{M}_{p(\cdot),q(\cdot),\lambda}^{loc}(\Rn)$, whose proof is the same of Theorem \ref{cay5}.

\begin{theorem}\label{cay57} Let $T_0$ be a sublinear operator satisfying condition  \eqref{BS1} and \eqref{sublR}. Let also   $1\leq p_{-}\leq p_{+}<\frac{q_*(t)}{\lambda}$, $1\leq q_{-}\leq q_{+}<\infty$, $p,q\in \mathcal{P}_{0,\i}(0,\i)$, $0 \leq \lambda < 1$ and $f\in \mathcal{M}_{p(\cdot),q(\cdot),\lambda}^{loc}(\Rn)$. Suppose that the conditions  \eqref{kuc22} and of Lemma \ref{cay0001} are satisfied. Then the operator $T_0$ is  bounded on the local variable Morrey-Lorentz spaces $\mathcal{M}_{p(\cdot),q(\cdot),\lambda}^{loc}(\Rn)$.
\end{theorem}

\section{Some Applications}

In this section we give some applications of our main results. We obtain the boundedness of Bochner-Riesz operator $B_r^\delta$, identity approximation $A_{\varepsilon}$ and the Marcinkiewicz operator $\mu_\Omega$ on the local variable Morrey-Lorentz spaces $\mathcal{M}_{p(\cdot),q(\cdot),\lambda}^{loc}(\Rn)$.

\,

\subsection{Bochner-Riesz operator }
Let $\delta>(n-1)/2$,
$B_r^\delta(f)^{\hat{}}(\xi)=(1-r^2|\xi|^2)_+^\delta\hat{f}(\xi)$
and $B_r^\delta(x)=r^{-n}B^\delta(x/r)$ for $r>0$. The maximal
Bochner-Riesz operator is defined by (see \cite{LiuLu}, \cite{LiuChen})
$$
B_{\delta,*}(f)(x)=\sup_{r>0}|B_{r}^\delta(f)(x)|.
$$


Let $H$ be the space $H=\{h:\|h\|=\sup_{t>0}|h(t)|<\infty\}$, then
it is clear that $B_{\delta,\ast}(f)(x)=\|B_{t}^\delta(f)(x)\|$.

By the condition on $B_r^\delta$ (see \cite{GarRub}), we have
\begin{align*}
	|B_r^\delta(x-y)|&\leq Cr^{-n}(1+|x-y|/r)^{-(\delta+(n+1)/2)}
	\\
	&=C\left(\frac{r}{r+|x-y|}\right)^{\delta-(n-1)/2}\frac{1}{(r+|x-y|)^n}
	\\
	&\leq |x-y|^{-n},
\end{align*}
and
\begin{equation*}
	B_{\delta,\ast}(f)(x)\leq C\int_{\Rn}\frac{|f(y)|}{|x-y|^n}dy.
\end{equation*}
Thus, the operator $B_{\delta,\ast}$ satisfies conditions \tcr{\eqref{BS1}} and \eqref{sublR}, then from Theorem \ref{cay57} we get
\begin{corollary} \label{ak13}  Let $1 < p_{-} \leq p_{+}< \infty , 1 \leq q_{-} \leq q_{+}< \infty $, $p, q\in \mathcal{P}_{0,\i}(0,\i)$, $0 \leq \lambda < 1$  and $f\in \mathcal{M}_{p(\cdot),q(\cdot),\lambda}^{loc}(\Rn)$.  Suppose that the conditions \eqref{kuc22} and of Lemma \ref{cay0001} are satisfied. Then  the Bochner-Riesz operator $B_r^\delta$ is bounded on the local variable Morrey-Lorentz spaces $\mathcal{M}_{p(\cdot),q(\cdot),\lambda}^{loc}(\Rn)$.
\end{corollary}

In the case $\lambda=0,$ from Theorem \ref{cay4} the boundedness of the operator $B_r^\delta$ in the variable Lorentz spaces $L_{p(\cdot),q(\cdot)}(\Rn)$ is obtained.


\subsection{Identity approximation }

It is known that the identity approximation (see \cite{St2})
\begin{align*}
A_{\varepsilon}f(x)=\frac{1}{\varepsilon^{n}}\int_{\Rn}a\left( \frac{x-y}{\varepsilon}\right) f(y)dy,
\end{align*}
where $\int_{\Rn}a(y)dy=1$ and $a(x)$ has a radial decreasing integrable majorant, are dominated by the maximal operator
\begin{align*}
|A_{\varepsilon}f(x)|\leq CMf(x), f\in L^{p}(\Rn), 1\leq p\leq \i,
\end{align*}
with an absolute constant $C>0$ not depending on $x$ and $\varepsilon$.

Since the maximal operator $M$ is bounded on the spaces ${\mathcal{M}_{p(\cdot),q(\cdot),\lambda}^{loc}}$, then we get the following.

\begin{corollary} \label{cay6} Let $1 < p_{-} \leq p_{+}< \infty , 1 \leq q_{-} \leq q_{+}< \infty $, $p,q\in \mathcal{P}_{0,\i}(0,\i)$,  $0 \leq \lambda < 1$ and $f\in \mathcal{M}_{p(\cdot),q(\cdot),\lambda}^{loc}(\Rn)$. Suppose that the conditions \eqref{kuc22} and of Lemma \ref{cay0001} are satisfied. Then the sublinear operator $\sup_{\varepsilon>0}|A_{\varepsilon}f(x)|$ is bounded on the local variable Morrey-Lorentz spaces $\mathcal{M}_{p(\cdot),q(\cdot),\lambda}^{loc}$.
\end{corollary}

In the case $\lambda=0,$ from Theorem \ref{cay4} we get the boundedness of the operator $\sup_{\varepsilon>0}|A_{\varepsilon}f(x)|$ in the variable Lorentz spaces $L_{p(\cdot),q(\cdot)}(\Rn)$ which is proved in \cite{EKS}.


\subsection{Marcinkiewicz operator}
In the case $n \geq 2$, we denote by $S^{n-1}=\{x\in\Rn:|x|=1\}$  the unit sphere in $\Rn$ equipped with the normalized Lebesgue measure $d\sigma$. Suppose that $\Omega$ satisfies the following conditions.

(i) $\Omega$ is the homogeneous function of degree zero on $\Rn\setminus\{0\}$, that is,
\begin{equation}\label{kuc05}
\Omega(t x)=\Omega(x),~~\text{for any}~~ t>0, \, x \in \Rn \setminus\{0\}.
\end{equation}

(ii) $\Omega$ has mean zero on $S^{n-1}$, that is,
\begin{equation}\label{kuc06}
\int_{S^{n-1}}\Omega(x^\prime)d\sigma(x^\prime)=0.
\end{equation}

(iii) $\Omega\in {\rm Lip}_\gamma(S^{n-1})$, $0<\gamma\leq1$, that is there exists a constant $C>0$ such that,
\begin{equation}\label{kuc07}
|\Omega(x^\prime)-\Omega(y^\prime)|\leq C |x^\prime-y^\prime|^\gamma~~\text{for any}~~ x^\prime,y^\prime\in S^{n-1}.
\end{equation}
where $x'=\frac{x}{|x|}$ for any $x\neq 0$.
In 1958, Stein \cite{Stein58} defined the Marcinkiewicz integral of higher dimension $\mu_\Omega$ satisfying the conditions \eqref{kuc05}-\eqref{kuc07} by
\begin{equation*}
\mu_{\Omega}(f)(x)=\left(\int_0^\infty|F_{\Omega,t}(f)(x)|^2\frac{dt}{t^3}\right)^{1/2},
\end{equation*}
where
$$
F_{\Omega,t}(f)(x)=\int_{|x-y|\leq t}\frac{\Omega(x-y)}{|x-y|^{n-1}}f(y)dy.
$$
The continuity of Marcinkiewicz operator $\mu_\Omega$ has been extensively studied in \cite{LuDingY, St2}.
Let $H$ be the space
$H=\{h:\|h\|=(\int_0^\infty|h(t)|^2dt/t^3)^{1/2}<\i\}$. Then, it is clear that $\mu_{\Omega}(f)(x)=\|F_{\Omega,t}(f)(x)\|$.

By Minkowski inequality and the above conditions on $\Omega$, we get
\begin{align*}
\mu_\Omega(f)(x) & = \left(\int_0^\infty \Big|\int_{\Rn}\frac{\Omega(x-y)}{|x-y|^{n-1}} \, \chi_{_{B(x,t)}}(y) \, f(y)dy\Big|^2\frac{dt}{t^3}\right)^{1/2}
\\
&\leq \int_{\Rn}\frac{|\Omega(x-y)|}{|x-y|^{n-1}}|f(y)|\left( \int_{0}^{\i} \, \chi_{_{B(x,t)}}(y) \,\frac{dt}{t^{3}} \right)^{\frac{1}{2}}dy
\\
&\leq \int_{\Rn}\frac{|\Omega(x-y)|}{|x-y|^{n-1}}|f(y)|\left( \int_{|x-y|}^{\i}\frac{dt}{t^{3}} \right)^{\frac{1}{2}}dy
\leq C \int_{\Rn}\frac{|f(y)|}{|x-y|^{n}}dy.
\end{align*}
Thus, the operator $\mu_\Omega$ satisfies conditions \tcr{\eqref{BS1}} and \eqref{sublR}, then from Theorem \ref{cay57} we get the following.
\begin{corollary} \label{kuc08} Let  $1 < p_{-} \leq p_{+}< \infty , 1 \leq q_{-} \leq q_{+}< \infty $, $p,q\in \mathcal{P}_{0,\i}(0,\i)$, $0 \leq \lambda < 1$ and $f\in \mathcal{M}_{p(\cdot),q(\cdot),\lambda}^{loc}(\Rn)$. Suppose that the conditions \eqref{kuc22} and of Lemma \ref{cay0001} are satisfied. Then the Marcinkiewicz operator $\mu_\Omega$ is  bounded on the local variable Morrey-Lorentz spaces $\mathcal{M}_{p(\cdot),q(\cdot),\lambda}^{loc}(\Rn)$.
\end{corollary}

In the case $\lambda=0,$ from Theorem \ref{cay5} we get the boundedness of the operator $\mu_\Omega$ in the variable Lorentz spaces $L_{p(\cdot),q(\cdot)}(\Rn)$.

\end{document}